\documentclass[12pt]{article}
\usepackage[a4paper, total={6in, 8in}, margin={1in}]{geometry}
\usepackage[utf8]{inputenc}

\usepackage{amsmath}
\usepackage{amsthm}
\usepackage{amsfonts}
\usepackage[hidelinks]{hyperref}
\usepackage{scrextend}

\theoremstyle{plain}
\newtheorem{theorem}{Theorem}
\newtheorem{corollary}[theorem]{Corollary}

\theoremstyle{definition}

\newtheorem{example}[theorem]{Example}

\theoremstyle{remark}
\newtheorem{remark}[theorem]{Remark}

\title{On an Identity by Bruckman and Good}
\author{Hongshen Chua}
\date{}

\begin{document}
\maketitle

\section{Introduction} \label{Sec: 1}

Define the \textit{generalized Lucas sequence}, $ (w_n)_{n \geq 0} = w_n(a, b; p, q) $ recursively by
\begin{align*}
	w_0 = a, \quad w_1 = b, \quad w_n = pw_{n - 1} - qw_{n - 2},
\end{align*}
where $ a, b, p $, and $ q $ can be complex numbers. Additionally, we define
\begin{align*}
	u_n (p, q) = w_n(0, 1; p, q), \quad v_n (p, q) = w_n(2, p; p, q),
\end{align*}
as the Lucas sequence of the first kind and second kind, respectively. The usual Fibonacci and Lucas numbers are given by $ F_n = u_n (1, -1) $ and $ L_n = v_n (1, -1) $, respectively. Denote $ \alpha $ and $ \beta $, with $ |\alpha| > |\beta| $, as the roots of the characteristic equation $ x^2 - px + q = 0 $. The discriminant of this equation is $ D = p^2 - 4q \neq 0 $. Notably, we have $ \alpha + \beta = p $, $ \alpha \beta = q $, and $ \alpha - \beta = \sqrt{D} $. The Binet formulas for $ w_n $, $ u_n $, and $ v_n $ are given by
\begin{align*}
	w_n = \frac{A\alpha^n - B\beta^n}{\alpha - \beta}, \quad u_n = \frac{\alpha^n - \beta^n}{\alpha - \beta}, \quad v_n = \alpha^n + \beta^n,
\end{align*}
where $ A = b - a\beta $ and $ B = b - a\alpha $.

Many papers have been devoted to finding closed forms for sums involving reciprocals of generalized Lucas sequences, with a particular focus on the Millin series. The series
\begin{align*}
	\sum_{i = 0}^\infty \frac{1}{F_{2^i}} = \frac{7 - \sqrt{5}}{2},
\end{align*}
was proposed by Miller \cite[p.\ 309]{Millin_74} as a problem in Fibonacci Quarterly. It is worth noting that this series appeared in earlier papers by Lucas \cite[Eq.\ (126)]{Lucas_65} and Brady \cite[Eq.\ (3)]{Brady_74}. In fact, they both gave a more general formula. For instance, Lucas demonstrated that
\begin{align*}
	\sum_{i = 0}^\infty \frac{q^{2^i r}}{u_{2^{i + 1} r}} = \frac{\beta^r}{u_r},
\end{align*}
where $ r \geq 1 $ is an integer. Other solutions include \cite{Bruckman_Good_76, Good_74, Hoggatt_Bicknell_76_2, Shannon_76}.

On top of that, numerous authors have also taken interest in the sum when the indices form an arithmetic progression, say
\begin{align*}
	\sum_{i = 1}^N \frac{q^{ri}}{w_{ri + s} w_{r(i + M) + s}},
\end{align*}
and its corresponding infinite version. Some examples include \cite{Adegoke_Frontczak_Goy_23, Hu_Sun_Liu_01, Melham_12}. Furthermore, Adegoke \cite{Adegoke_18} extended the sum to cases where the denominator contains two or more terms.

Traditionally, one can derive these results by applying telescoping sums to recurrence relations of generalized Lucas sequences. An alternative was provided by Bruckman and Good \cite{Bruckman_Good_76}, via an identity that can be traced back to de Morgan, i.e.,
\begin{align*}
	\frac{x_{k + 1} y_k - x_k y_{k + 1}}{(x_k - y_k)(x_{k + 1} - y_{k + 1})} = \frac{y_k}{x_k - y_k} - \frac{y_{k + 1}}{x_{k + 1} - y_{k + 1}}.
\end{align*}
However, Bruckman and Good only gave sums involving Fibonacci and Lucas numbers. Further extensions to Lucas sequences were added by Farhi \cite{Farhi_19}. Finally, we have to mention that de Morgan's identity is a special case of a more general identity by Duverney and Shiokawa \cite[Cor.\ 4.1]{Duverney_Shiokawa_76}. It could be worthwhile to use this identity to discover even more extensive series involving Lucas sequences.

In this paper, we consider two variations on the identity by Bruckman and Good, enabling us to establish numerous intriguing sums and rediscover some previously known ones. We will make repeated use of the following telescoping summation identities. For any positive integers $ N $ and $ M $,
\begin{align} \label{Eq: Telescoping One}
	\sum_{i = 1}^N (x_i - x_{i + M}) = \sum_{i = 1}^M (x_i - x_{i + N}),
\end{align}
and 
\begin{align} \label{Eq: Telescoping Two}
	\sum_{i = 1}^{2N} (\pm 1)^i (x_i - x_{i + 2M}) = \sum_{i = 1}^{2M} (\pm 1)^i (x_i - x_{i + 2N}).
\end{align}
Note that the usual telescoping sum is obtained by setting $ M = 1 $ in Eq.\ (\ref{Eq: Telescoping One}). 

\section{New reciprocal series of Lucas sequences} \label{Sec: 2}

Let $ A $ and $ B $ be arbitrary complex numbers. The identity provided by Bruckman and Good \cite{Bruckman_Good_76} can be generalized for any sequences $(x_n)$ and $(y_n)$ as follows:
\begin{align*}
	\frac{A(x_{i + j} y_i - x_i y_{i + j})}{(Ax_i - By_i)(Ax_{i + j} - By_{i + j})} = \frac{y_i}{Ax_i - By_i} - \frac{y_{i + j}}{Ax_{i + j} - By_{i + j}}.
\end{align*}
Since the RHS is telescoping, upon applications of Eq.\ (\ref{Eq: Telescoping One}) and Eq.\ (\ref{Eq: Telescoping Two}), we obtain
\begin{align*}
	\sum_{i = 1}^N \frac{A(x_{i + M} y_i - x_i y_{i + M})}{(Ax_i - By_i)(Ax_{i + M} - By_{i + M})} &= \sum_{i = 1}^M \frac{A(x_{i + N} y_i - x_i y_{i + N})}{(Ax_i - By_i)(Ax_{i + N} - By_{i + N})}, \\
	\sum_{i = 1}^{2N} \frac{(\pm 1)^i A(x_{i + 2M} y_i - x_i y_{i + 2M})}{(Ax_i - By_i)(Ax_{i + 2M} - By_{i + 2M})} &= \sum_{i = 1}^{2M} \frac{(\pm 1)^i  A(x_{i + 2N} y_i - x_i y_{i + 2N})}{(Ax_i - By_i)(Ax_{i + 2N} - By_{i + 2N})}.
\end{align*}
In particular, for any integer function $ f $,
\begin{align*}
	\sum_{i = 1}^N \frac{A(x^{f(i + M)} y^{f(i)} - x^{f(i)} y^{f(i + M)})}{(Ax^{f(i)} - By^{f(i)})(Ax^{f(i + M)} - By^{f(i + M)})} &= \sum_{i = 1}^M \frac{A(x^{f(i + N)} y^{f(i)} - x^{f(i)} y^{f(i + N)})}{(Ax^{f(i)} - By^{f(i)})(Ax^{f(i + N)} - By^{f(i + N)})}, \\
	\sum_{i = 1}^{2N} \frac{(\pm 1)^i A(x^{f(i + 2M)} y^{f(i)} - x^{f(i)} y^{f(i + 2M)})}{(Ax^{f(i)} - By^{f(i)})(Ax^{f(i + 2M)} - By^{f(i + 2M)})} &= \sum_{i = 1}^{2M} \frac{(\pm 1)^i A(x^{f(i + 2N)} y^{f(i)} - x^{f(i)} y^{f(i + 2N)})}{(Ax^{f(i)} - By^{f(i)})(Ax^{f(i + 2N)} - By^{f(i + 2N)})}.
\end{align*}
By setting $ x = \alpha $ and $ y = \beta $, we obtain the following theorems:
\begin{theorem} \label{Th: 1} For integer function $ f $ and positive integers $ M $ and $ N $,
\begin{align*}
	\sum_{i = 1}^N \frac{q^{f(i)} u_{f(i + M) - f(i)}}{w_{f(i)} w_{f(i + M)}} &= \sum_{i = 1}^M \frac{q^{f(i)} u_{f(i + N) - f(i)}}{w_{f(i)} w_{f(i + N)}}. 
\end{align*}
\end{theorem}

\begin{theorem} \label{Th: 2} For integer function $ f $ and positive integers $ M $ and $ N $,
\begin{align*}
	\sum_{i = 1}^{2N} \frac{(\pm 1)^i q^{f(i)} u_{f(i + 2M) - f(i)}}{w_{f(i)} w_{f(i + 2M)}} &= \sum_{i = 1}^{2M} \frac{(\pm 1)^i q^{f(i)} u_{f(i + 2N) - f(i)}}{w_{f(i)} w_{f(i + 2N)}}. 
\end{align*}
\end{theorem}

In particular, if $ M = 1 $, then  
\begin{align*}
	\sum_{i = 1}^N \frac{Aq^{f(i)} u_{f(i + 1) - f(i)}}{w_{f(i)} w_{f(i + 1)}} &= \frac{Aq^{f(1)} u_{f(N + 1) - f(1)}}{w_{f(1)} w_{f(N + 1)}} = \frac{\beta^{f(1)}}{w_{f(1)}} - \frac{\beta^{f(N + 1)}}{w_{f(N + 1)}}, 
\end{align*}
and
\begin{align*}
	\sum_{i = 1}^{2N} \frac{(\pm 1)^i Aq^{s(i)} u_{f(i + 2) - f(i)}}{w_{f(i)} w_{f(i + 2)}} &= \pm \frac{Aq^{f(1)} u_{f(2N + 1) - f(1)}}{w_{f(1)} w_{f(2N + 1)}} + \frac{Aq^{f(2)} u_{f(2N + 2) - f(2)}}{w_{f(2)} w_{f(2N + 2)}} \\
	&= \pm \frac{\beta^{f(1)}}{w_{f(1)}} \mp \frac{\beta^{f(2N + 1)}}{w_{f(2N + 1)}} + \frac{\beta^{f(2)}}{w_{f(2)}} - \frac{\beta^{f(2N + 2)}}{w_{f(2N + 2)}}.
\end{align*}
Assuming further that $ f(i) \rightarrow \infty $ as $ i \rightarrow \infty $. By letting $ N \rightarrow \infty $, each term containing $ N $ on the RHS vanishes, so we get the following infinite versions of the previous theorems.
\begin{theorem} \label{Th: 3} Suppose that $ f(i) \rightarrow \infty $ as $ i \rightarrow \infty $. Then
\begin{align*}
	\sum_{i = 1}^\infty \frac{q^{f(i)} u_{f(i + 1) - f(i)}}{w_{f(i)} w_{f(i + 1)}} &= \frac{\beta^{f(1)}}{A w_{f(1)}}. 
\end{align*}
\end{theorem}

\begin{remark} This was also derived by Hu et al.\ \cite[Th.\ 2]{Hu_Sun_Liu_01}. For increments greater than $ 1 $, refer to Farhi's paper \cite[Th. 2]{Farhi_19}.
\end{remark}

\begin{theorem} \label{Th: 5} Suppose that $ f(i) \rightarrow \infty $ as $ i \rightarrow \infty $. Then
\begin{align*}
	\sum_{i = 1}^\infty \frac{(\pm 1)^i q^{f(i)} u_{f(i + 2) - f(i)}}{w_{f(i)} w_{f(i + 2)}} &= \frac{\beta^{f(2)}}{Aw_{f(2)}} \pm \frac{\beta^{f(1)}}{Aw_{f(1)}}. 
\end{align*}
\end{theorem}

It is immediate to deduce from Theorem \ref{Th: 5} the next result.
\begin{theorem} \label{Th: 6} Suppose that $ f(i) \rightarrow \infty $ as $ i \rightarrow \infty $. Then
\begin{align*}
	\sum_{\substack{i = 1 \\ i \text{ even}}}^\infty \frac{q^{f(i)} u_{f(i + 2) - f(i)}}{w_{f(i)} w_{f(i + 2)}} = \frac{\beta^{f(2)}}{Aw_{f(2)}}, \\
	\sum_{\substack{i = 1 \\ i \text{ odd}}}^\infty \frac{q^{f(i)} u_{f(i + 2) - f(i)}}{w_{f(i)} w_{f(i + 2)}} = \frac{\beta^{f(1)}}{Aw_{f(1)}}. 
\end{align*}
\end{theorem}

\section{Applications, Part I} \label{Sec: 3}

\subsection{When $ f(i) = ri + s $ with integers $ r, s $ and $ r > 0 $}

Using Theorems \ref{Th: 1} to \ref{Th: 5}, we obtain the next corollaries as special cases.
\begin{corollary} \label{Cor: 7} Let $ r $, $ s $, $ M $, and $ N $ be integers with $ r $, $ M $, and $ N > 0 $. Then
\begin{align*}
	u_{rM} \sum_{i = 1}^N \frac{q^{ri}}{w_{ri + s} w_{r(i + M) + s}} &= u_{rN} \sum_{i = 1}^M \frac{q^{ri}}{w_{ri + s} w_{r(i + N) + s}}, \\
	u_{2rM} \sum_{i = 1}^{2N} \frac{(\pm 1)^i q^{ri}}{w_{ri + s} w_{r(i + 2M) + s}} &= u_{2rN} \sum_{i = 1}^{2M} \frac{(\pm 1)^i q^{ri}}{w_{ri + s} w_{r(i + 2N) + s}}. 
\end{align*}
\end{corollary}

\begin{corollary} \label{Cor: 8} Let $ r $ and $ s $ be integers with $ r > 0 $. Then
\begin{align*}
	\sum_{i = 1}^\infty \frac{q^{ri + s}}{w_{ri + s} w_{r(i + 1) + s}} &= \frac{\beta^{r + s}}{A u_r w_{r + s}}, \\
	\sum_{i = 1}^\infty \frac{(\pm 1)^i q^{ri + s}}{w_{ri + s} w_{r(i + 2) + s}} &= \frac{\beta^{2r + s}}{Au_{2r} w_{2r + s}} \pm \frac{\beta^{r + s}}{Au_{2r} w_{r + s}}.
\end{align*}
\end{corollary}

\begin{remark} These corollaries include the results of Adegoke et al.\ \cite{Adegoke_Frontczak_Goy_23}. For example, if we carry out the transformations $ M \rightarrow 2M $ and $ s \rightarrow s - rM $, then our first sum in Corollary \ref{Cor: 7} becomes
\begin{align*}
	u_{2rM} \sum_{i = 1}^N \frac{q^{ri}}{w_{r(i - M) + s} w_{r(i + M) + s}} &= u_{rN} \sum_{i = 1}^{2M} \frac{q^{ri}}{w_{r(i - M) + s} w_{r(i + N - M) + s}}.
\end{align*}
This is equivalent to Theorem $ 1 $ of the aforementioned paper. The rest can be deduced analogously. 
\end{remark}

\begin{example} \label{Eg: 10} To apply to the case of Fibonacci and Lucas numbers, we set $ p = 1 $ and $ q = -1 $, then $ \alpha = \varphi = (1 + \sqrt{5})/2 $, $ \beta = 1 - \varphi = (1 - \sqrt{5})/2 $, and $ D = 5 $. Corollary \ref{Cor: 8} gives
\begin{align*}
	\sum_{i = 1}^\infty \frac{(-1)^{ri + s}}{F_{ri + s} F_{r(i + 1) + s}} &= \frac{(1 - \varphi)^{r + s}}{F_r F_{r + s}}, \\
	\sum_{i = 1}^\infty \frac{(-1)^{ri + s}}{L_{ri + s} L_{r(i + 1) + s}} &= \frac{(1 - \varphi)^{r + s}}{\sqrt{5} L_r L_{r + s}}, \\
	\sum_{i = 1}^\infty \frac{(\pm 1)^i (-1)^{ri + s}}{F_{ri + s} F_{r(i + 2) + s}} &= \frac{(1 - \varphi)^{2r + s}}{F_{2r} F_{2r + s}} \pm \frac{(1 - \varphi)^{r + s}}{F_{2r} F_{r + s}}, \\
	\sum_{i = 1}^\infty \frac{(\pm 1)^i (-1)^{ri + s}}{L_{ri + s} L_{r(i + 2) + s}} &= \frac{(1 - \varphi)^{2r + s}}{\sqrt{5} F_{2r} L_{2r + s}} \pm \frac{(1 - \varphi)^{r + s}}{\sqrt{5} F_{2r} L_{r + s}}. 
\end{align*}

\begin{remark} The first two sums, when both $ r $ and $ s $ are even, can be found in Popov's paper \cite[p.\ 263]{Popov_84}. The general case for $ r $ being even and $ s \neq 0 $ was proved by Melham \cite[Th.\ 2, 3]{Melham_12}. 
\end{remark}
\end{example} 

\begin{example} \label{Eg: 11} Set $ p = e + 1/e $, $ q = 1 $, $ A = 1 $, and $ B = \pm 1 $, then $ \alpha = e $ and $ \beta = 1 /e $. Corollary \ref{Cor: 8} gives
\begin{align*}
	\sum_{i = 1}^\infty \frac{\sinh r }{\sinh (ri + s) \sinh (r(i + 1) + s)} &= \frac{1}{e^{r + s} \sinh (r + s)}, \\
	\sum_{i = 1}^\infty \frac{\sinh r }{\cosh (ri + s) \cosh (r(i + 1) + s)} &= \frac{1}{e^{r + s} \cosh (r + s)}, \\
	\sum_{i = 1}^\infty \frac{(\pm 1)^i \sinh 2r}{\sinh (ri + s) \sinh (r(i + 2) + s)} &= \frac{1}{e^{2r + s} \sinh (2r + s)} \pm \frac{1}{e^{r + s} \sinh (r + s)}, \\
	\sum_{i = 1}^\infty \frac{(\pm 1)^i \sinh 2r}{\cosh (ri + s) \cosh (r(i + 2) + s)} &= \frac{1}{e^{2r + s} \cosh (2r + s)} \pm \frac{1}{e^{r + s} \cosh (r + s)}. 
\end{align*}
\end{example} 	

\subsection{When $ f(i) = 2^i r $ with integer $ r > 0 $} 

Using Theorems \ref{Th: 1} to \ref{Th: 5}, we have the next results.
\begin{corollary} \label{Cor: 13} Let $ r $, $ s $, $ M $, and $ N $ be integers with $ r $, $ M $, and $ N > 0 $. Then
\begin{align*}
	\sum_{i = 1}^N \frac{q^{2^i r} u_{2^{i + M} r - 2^i r}}{w_{2^i r} w_{2^{i + M} r}} &= \sum_{i = 1}^M \frac{q^{2^i r} u_{2^{i + N} r - 2^i r}}{w_{2^i r} w_{2^{i + N} r}}, \\
	\sum_{i = 1}^{2N} \frac{(\pm 1)^i q^{2^i r} u_{2^{i + 2M} r - 2^i r}}{w_{2^i r} w_{2^{i + 2M} r}} &= \sum_{i = 1}^{2M} \frac{(\pm 1)^i q^{2^i r} u_{2^{i + 2N} r - 2^i r}}{w_{2^i r} w_{2^{i + 2N} r}}. 
\end{align*}
\end{corollary}

\begin{corollary} \label{Cor: 14} Let $ r > 0 $ be an integer. Then
\begin{align*}
	\sum_{i = 1}^\infty \frac{q^{2^i r} u_{2^i r}}{w_{2^i r} w_{2^{i + 1} r}} &= \frac{\beta^{2r}}{Aw_{2r}}, \\
	\sum_{i = 1}^\infty \frac{(\pm 1)^i q^{2^i r} u_{3 \cdot 2^i r}}{w_{2^i r} w_{2^{i + 2} r}} &= \frac{\beta^{4r}}{Aw_{4r}} \pm \frac{\beta^{2r}}{Aw_{2r}}.
\end{align*}
\end{corollary}

\begin{example} \label{Eg: 15} To apply to the case of Fibonacci and Lucas numbers, we set $ p = 1 $ and $ q = -1 $, then $ \alpha = \varphi $, $ \beta = 1 - \varphi $, and $ D = 5 $. Corollary \ref{Cor: 14} gives
\begin{align*}
	\sum_{i = 1}^\infty \frac{1}{F_{2^{i + 1} r}} &= \frac{(1 - \varphi)^{2r}}{F_{2r}}, \\
	\sum_{i = 1}^\infty \frac{F_{2^i r}}{L_{2^i r} L_{2^{i + 1} r}} &= \frac{(1 - \varphi)^{2r}}{\sqrt{5} L_{2r}}, \\
	\sum_{i = 1}^\infty \frac{(\pm 1)^i F_{3 \cdot 2^i r}}{F_{2^i r} F_{2^{i + 2} r}} &= \frac{(1 - \varphi)^{4r}}{F_{4r}} \pm \frac{(1 - \varphi)^{2r}}{F_{2r}}, \\
	\sum_{i = 1}^\infty \frac{(\pm 1)^i F_{3 \cdot 2^i r}}{L_{2^i r} L_{2^{i + 2} r}} &= \frac{(1 - \varphi)^{4r}}{\sqrt{5} L_{4r}} \pm \frac{(1 - \varphi)^{2r}}{\sqrt{5} L_{2r}}. 
\end{align*}
\end{example}

\begin{example} \label{Eg: 16} Set $ p = e + 1/e $, $ q = 1 $, $ A = 1 $, and $ B = \pm 1 $, then $ \alpha = e $ and $ \beta = 1 /e $. Corollary \ref{Cor: 14} gives
\begin{align*}
	\sum_{i = 1}^\infty \frac{1}{\sinh 2^{i + 1} r} &= \frac{1}{e^{2r} \sinh 2r}, \\
	\sum_{i = 1}^\infty \frac{\tanh 2^i r}{\cosh 2^{i + 1} r} &= \frac{1}{e^{2r} \cosh 2r}, \\
	\sum_{i = 1}^\infty \frac{(\pm 1)^i \sinh (3 \cdot 2^i r)}{\sinh 2^i r \sinh 2^{i + 2} r} &= \frac{1}{e^{4r} \sinh 4r} \pm \frac{1}{e^{2r} \sinh 2r}, \\
	\sum_{i = 1}^\infty \frac{(\pm 1)^i \sinh (3 \cdot 2^i r)}{\cosh 2^i r \cosh 2^{i + 2} r} &= \frac{1}{e^{4r} \cosh 4r} \pm \frac{1}{e^{2r} \cosh 2r}.
\end{align*}
\end{example} 

\begin{remark} The first sum was also discovered by Gould \cite[Eq.\ 24]{Gould_77}.	
\end{remark}

\subsection{When $ f(i) = 3^i r $ with integer $ r > 0 $} 

The application of Theorems \ref{Th: 1} to \ref{Th: 5} yields the next two corollaries.
\begin{corollary} \label{Cor: 18} Let $ r $, $ s $, $ M $, and $ N $ be integers with $ r $, $ M $, and $ N > 0 $. Then
\begin{align*}
	\sum_{i = 1}^N \frac{q^{3^i r} u_{3^{i + M} r - 3^i r}}{w_{3^i r} w_{3^{i + M} r}} &= \sum_{i = 1}^M \frac{q^{3^i r} u_{3^{i + N} r - 3^i r}}{w_{3^i r} w_{3^{i + N} r}}, \\
	\sum_{i = 1}^{2N} \frac{(\pm 1)^i q^{3^i r} u_{3^{i + 2M} r - 3^i r}}{w_{3^i r} w_{3^{i + 2M} r}} &= \sum_{i = 1}^{2M} \frac{(\pm 1)^i q^{3^i r} u_{3^{i + 2N} r - 3^i r}}{w_{3^i r} w_{3^{i + 2N} r}}.
\end{align*}
\end{corollary}

\begin{corollary} \label{Cor: 19} Let $ r > 0 $ be an integer. Then
\begin{align*}
	\sum_{i = 1}^\infty \frac{q^{3^i r} u_{2 \cdot 3^i r}}{w_{3^i r} w_{3^{i + 1} r}} &= \frac{\beta^{3r}}{Aw_{3r}}, \\
	\sum_{i = 1}^\infty \frac{(\pm 1)^i q^{3^i r} u_{8 \cdot 3^i r}}{w_{3^i r} w_{3^{i + 2} r}} &= \frac{\beta^{9r}}{Aw_{9r}} \pm \frac{\beta^{3r}}{Aw_{3r}}.
\end{align*}
\end{corollary}

\begin{example} \label{Eg: 20} To apply to the case of Fibonacci and Lucas numbers, we set $ p = 1 $ and $ q = -1 $, then $ \alpha = \varphi $, $ \beta = 1 - \varphi $, and $ D = 5 $. Corollary \ref{Cor: 19} gives
\begin{align*}
	\sum_{i = 1}^\infty \frac{(-1)^r L_{3^i r}}{F_{3^{i + 1} r}} &= \frac{(1 - \varphi)^{3r}}{F_{3r}}, \\
	\sum_{i = 1}^\infty \frac{(-1)^r F_{3^i r}}{L_{3^{i + 1} r}} &= \frac{(1 - \varphi)^{3r}}{\sqrt{5} L_{3r}}, \\
	\sum_{i = 1}^\infty \frac{(\pm 1)^i (-1)^r F_{8 \cdot 3^i r}}{F_{3^i r} F_{3^{i + 2} r}} &= \frac{(1 - \varphi)^{9r}}{\sqrt{5} F_{9r}} \pm \frac{(1 - \varphi)^{3r}}{\sqrt{5} F_{3r}}, \\
	\sum_{i = 1}^\infty \frac{(\pm 1)^i (-1)^r F_{8 \cdot 3^i r}}{L_{3^i r} L_{3^{i + 2} r}} &= \frac{(1 - \varphi)^{9r}}{\sqrt{5} L_{9r}} \pm \frac{(1 - \varphi)^{3r}}{\sqrt{5} L_{3r}}. 
\end{align*}
\end{example} 

\begin{remark} The first two sums appeared in the paper by Bruckman and Good \cite[Eq.\ (11)]{Bruckman_Good_76}. Shar \cite[p.\ 10]{Shar_12} also discovered a similar series.
\end{remark}

\begin{example} \label{Eg: 22} Set $ p = e + 1/e $, $ q = 1 $, $ A = 1 $, and $ B = \pm 1 $, then $ \alpha = e $ and $ \beta = 1 /e $. Corollary \ref{Cor: 19} gives
\begin{align*}
	\sum_{i = 1}^\infty \frac{\cosh 3^i r}{\sinh 3^{i + 1} r} &= \frac{1}{2e^{3r} \sinh 3r}, \\
	\sum_{i = 1}^\infty \frac{\sinh 3^i r}{\cosh 3^{i + 1} r} &= \frac{1}{2e^{3r} \cosh 3r}, \\
	\sum_{i = 1}^\infty \frac{(\pm 1)^i \sinh (8 \cdot 3^i r)}{\sinh 3^i r \sinh 3^{i + 2} r} &= \frac{1}{e^{9r} \sinh 9r} \pm \frac{1}{e^{3r} \sinh 3r}, \\
	\sum_{i = 1}^\infty \frac{(\pm 1)^i \sinh (8 \cdot 3^i r)}{\cosh 3^i r \cosh 3^{i + 2} r} &= \frac{1}{e^{9r} \cosh 9r} \pm \frac{1}{e^{3r} \cosh 3r}.
\end{align*}
\end{example} 

\section{More reciprocal series of Lucas sequences} \label{Sec: 4}

In this section, we will only deal with the Lucas sequences, $ (u_n) $ and $ (v_n) $, instead of the generalized Lucas sequences, $ (w_n) $. Given any sequences $ (x_n) $ and $ (y_n) $, one checks that
\begin{align*}
	\frac{\pm 4(x_{i + j} y_i - x_i y_{i + j})(x_i x_{i + j} - y_i y_{i + j})}{(x_i \mp y_i)^2 (x_{i + j} \mp y_{i + j})^2} = \frac{(x_i \pm y_i)^2}{(x_i \mp y_i)^2} - \frac{(x_{i + j} \pm y_{i + j})^2}{(x_{i + j} \mp y_{i + j})^2}.
\end{align*}
Using the same procedure as in \S \ref{Sec: 2}, we deduce the next two theorems.
\begin{theorem} \label{Th: 23} For integer function $ f $ and positive integers $ N $ and $ M $,
\begin{align*}
	\sum_{i = 1}^N \frac{q^{f(i)} u_{f(i + M) - f(i)} u_{f(i + M) + f(i)}}{u_{f(i)}^2 u_{f(i + M)}^2} &= \sum_{i = 1}^M \frac{q^{f(i)} u_{f(i + N) - f(i)} u_{f(i + N) + f(i)}}{u_{f(i)}^2 u_{f(i + N)}^2}, \\
	\sum_{i = 1}^N \frac{q^{f(i)} u_{f(i + M) - f(i)} u_{f(i + M) + f(i)}}{v_{f(i)}^2 v_{f(i + M)}^2} &= \sum_{i = 1}^M \frac{q^{f(i)} u_{f(i + N) - f(i)} u_{f(i + N) + f(i)}}{v_{f(i)}^2 v_{f(i + N)}^2}. 
\end{align*}
\end{theorem}

\begin{theorem} \label{Th: 24} For integer function $ f $ and positive integers $ N $ and $ M $,
\begin{align*}
	\sum_{i = 1}^{2N} \frac{(\pm 1)^i q^{f(i)} u_{f(i + 2M) - f(i)} u_{f(i + 2M) + f(i)}}{u_{f(i)}^2 u_{f(i + 2M)}^2} &= \sum_{i = 1}^{2M} \frac{(\pm 1)^i q^{f(i)} u_{f(i + 2N) - f(i)} u_{f(i + 2N) + f(i)}}{u_{f(i)}^2 u_{f(i + 2N)}^2}, \\
	\sum_{i = 1}^{2N} \frac{(\pm 1)^i q^{f(i)} u_{f(i + 2M) - f(i)} u_{f(i + 2M) + f(i)}}{v_{f(i)}^2 v_{f(i + 2M)}^2} &= \sum_{i = 1}^{2M} \frac{(\pm 1)^i q^{f(i)} u_{f(i + 2N) - f(i)} u_{f(i + 2N) + f(i)}}{v_{f(i)}^2 v_{f(i + 2N)}^2}. 
\end{align*}
\end{theorem}
In Theorems \ref{Th: 25} and \ref{Th: 26}, we assume further that $ f(i) \rightarrow \infty $ as $ i \rightarrow \infty $. When $ M = 1 $, using the fact that $ v_N^2/u_N^2 \rightarrow D $ as $ N \rightarrow \infty $, Theorems \ref{Th: 23} and \ref{Th: 24} imply that
\begin{theorem} \label{Th: 25} Suppose that $ f(i) \rightarrow \infty $ as $ i \rightarrow \infty $. Then
\begin{align*}
	\sum_{i = 1}^\infty \frac{4q^{f(i)} u_{f(i + 1) - f(i)} u_{f(i + 1) + f(i)}}{u_{f(i)}^2 u_{f(i + 1)}^2} &= \frac{v_{f(1)}^2}{u_{f(1)}^2} - D, \\
	\sum_{i = 1}^\infty \frac{4q^{f(i)} u_{f(i + 1) - f(i)} u_{f(i + 1) + f(i)}}{v_{f(i)}^2 v_{f(i + 1)}^2} &= \frac{1}{D} - \frac{u_{f(1)}^2}{v_{f(1)}^2}.
\end{align*}
\end{theorem} 

\begin{theorem} \label{Th: 26} Suppose that $ f(i) \rightarrow \infty $ as $ i \rightarrow \infty $. Then
\begin{align*}
	\sum_{i = 1}^\infty \frac{(\pm 1)^i 4 q^{f(i)} u_{f(i + 2) - f(i)} u_{f(i + 2) + f(i)}}{u_{f(i)}^2 u_{f(i + 2)}^2} &= \frac{v_{f(2)}^2}{u_{f(2)}^2} \pm \frac{v_{f(1)}^2}{u_{f(1)}^2} - D (1 \pm 1), \\
	\sum_{i = 1}^\infty \frac{(\pm 1)^i 4 q^{f(i)} u_{f(i + 2) - f(i)} u_{f(i + 2) + f(i)}}{v_{f(i)}^2 v_{f(i + 2)}^2} &= -\frac{u_{f(2)}^2}{v_{f(2)}^2} \mp \frac{u_{f(1)}^2}{v_{f(1)}^2} + \frac{1}{D} (1 \pm 1).
\end{align*}
\end{theorem}

\section{Applications, Part II} \label{Sec: 5}

\subsection{When $ f(i) = ri + s $ with integers $ r, s $ and $ r > 0 $}

The next results are direct consequences of Theorems \ref{Th: 23} to \ref{Th: 26}.
\begin{corollary} \label{Cor: 27} Let $ r $, $ s $, $ M $, and $ N $ be integers with $ r $, $ M $, and $ N > 0 $. Then
\begin{align*}
	u_{rM} \sum_{i = 1}^N \frac{q^{ri + s} u_{r(2i + M) + 2s}}{u_{ri + s}^2 u_{r(i + M) + s}^2} &= u_{rN} \sum_{i = 1}^M \frac{q^{ri + s} u_{r(2i + N) + 2s}}{u_{ri + s}^2 u_{r(i + N) + s}^2}, \\
	u_{rM} \sum_{i = 1}^N \frac{q^{ri + s} u_{r(2i + M) + 2s}}{v_{ri + s}^2 v_{r(i + M) + s}^2} &= u_{rN} \sum_{i = 1}^M \frac{q^{ri + s} u_{r(2i + N) + 2s}}{v_{ri + s}^2 v_{r(i + N) + s}^2}, \\
	u_{2rM} \sum_{i = 1}^{2N} \frac{(\pm 1)^i q^{ri + s} u_{2r(i + M) + 2s}}{u_{ri + s}^2 u_{r(i + 2M) + s}^2} &= u_{2rN} \sum_{i = 1}^{2M} \frac{(\pm 1)^i q^{ri + s} u_{2r(i + N) + 2s}}{u_{ri + s}^2 u_{r(i + 2N) + s}^2}, \\
	u_{2rM} \sum_{i = 1}^{2N} \frac{(\pm 1)^i q^{ri + s} u_{2r(i + M) + 2s}}{v_{ri + s}^2 v_{r(i + 2M) + s}^2} &= u_{2rN} \sum_{i = 1}^{2M} \frac{(\pm 1)^i q^{ri + s} u_{2r(i + N) + 2s}}{v_{ri + s}^2 v_{r(i + 2N) + s}^2}.
\end{align*}
\end{corollary}

\begin{corollary} \label{Cor: 28} Let $ r $ and $ s $ be integers with $ r > 0 $. Then
\begin{align*}
	\sum_{i = 1}^\infty \frac{4q^{ri + s} u_r u_{r(2i + 1) + 2s}}{u_{ri + s}^2 u_{r(i + 1) + s}^2} &= \frac{v_{r + s}^2}{u_{r + s}^2} - D, \\
	\sum_{i = 1}^\infty \frac{4q^{ri + s} u_r u_{r(2i + 1) + 2s}}{v_{ri + s}^2 v_{r(i + 1) + s}^2} &= \frac{1}{D} - \frac{u_{r + s}^2}{v_{r + s}^2}, \\
	\sum_{i = 1}^\infty \frac{(\pm 1)^i 4q^{ri + s} u_{2r} u_{2r(i + 1) + 2s}}{u_{ri + s}^2 u_{r(i + 2) + s}^2} &= \frac{v_{2r + s}^2}{u_{2r + s}^2} \pm \frac{v_{r + s}^2}{u_{r + s}^2} - D(1 \pm 1 ), \\
	\sum_{i = 1}^\infty \frac{(\pm 1)^i 4q^{ri + s} u_{2r} u_{2r(i + 1) + 2s}}{v_{ri + s}^2 v_{r(i + 2) + s}^2} &= -\frac{u_{2r + s}^2}{v_{2r + s}^2} \mp \frac{u_{r + s}^2}{v_{r + s}^2} + \frac{1}{D} (1 \pm 1).
\end{align*}
\end{corollary}

\begin{example} \label{Eg: 29} To apply to the case of Fibonacci and Lucas numbers, we set $ p = 1 $ and $ q = -1 $, then $ \alpha = \varphi $, $ \beta = 1 - \varphi $, and $ D = 5 $. Corollary \ref{Cor: 28} gives
\begin{align*}
	\sum_{i = 1}^\infty \frac{(-1)^{ri + s} F_{r(2i + 1) + 2s}}{F_{ri + s}^2 F_{r(i + 1) + s}^2} &= \frac{L_{r + s}^2}{4F_r F_{r + s}^2} - \frac{5}{4F_r}, \\
	\sum_{i = 1}^\infty \frac{(-1)^{ri + s} F_{r(2i + 1) + 2s}}{L_{ri + s}^2 L_{r(i + 1) + s}^2} &= \frac{1}{20F_r} - \frac{F_{r + s}^2}{4F_r L_{r + s}^2}, \\
	\sum_{i = 1}^\infty \frac{(\pm 1)^i (-1)^{ri + s} F_{2r(i + 1) + 2s}}{F_{ri + s}^2 F_{r(i + 2) + s}^2} &= \frac{L_{2r + s}^2}{4F_r F_{2r + s}^2} \pm \frac{L_{r + s}^2}{4F_r F_{r + s}^2} - \frac{5}{4F_r} (1 \pm 1), \\
	\sum_{i = 1}^\infty \frac{(\pm 1)^i (-1)^{ri + s} F_{2r(i + 1) + 2s}}{L_{ri + s}^2 L_{r(i + 2) + s}^2} &= -\frac{F_{2r + s}^2}{4F_r L_{2r + s}^2} \mp \frac{F_{r + s}^2}{4F_r L_{r + s}^2} + \frac{1}{20F_r}(1 \pm 1).
\end{align*}
\end{example}

\begin{remark} The first two sums were a problem in Fibonacci Quarterly, proposed and solved by Gauthier \cite{Gauthier_08, Gauthier_11}. However, Gauthier made an error in the sign of the first sum. We have also found several papers discussing similar results, but to our knowledge, none of them provided the general form as in Corollary \ref{Cor: 28}. The list of sources (not exhaustive) includes \cite{Adegoke_18, Adegoke_Frontczak_Goy_21, Brousseau_69, Bruckman_11, Hendel_12, Koshy_14, Plaza_Falcon_11}. 
\end{remark}

\subsection{When $ f(i) = 2^i r $ with integer $ r > 0 $} 

Applying Theorems \ref{Th: 23} to \ref{Th: 26}, we get the next two corollaries.
\begin{corollary} \label{Cor: 31} Let $ r $, $ s $, $ M $, and $ N $ be integers with $ r $, $ M $, and $ N > 0 $. Then
\begin{align*}
	\sum_{i = 1}^N \frac{q^{2^i r} u_{2^{i + M} r - 2^i r} u_{2^{i + M} r + 2^i r}}{u_{2^i r}^2 u_{2^{i + M} r}^2} &= \sum_{i = 1}^M \frac{q^{2^i r} u_{2^{i + N} r - 2^i r} u_{2^{i + N} r + 2^i r}}{u_{2^i r}^2 u_{2^{i + N} r}^2}, \\
	\sum_{i = 1}^N \frac{q^{2^i r} u_{2^{i + M} r - 2^i r} u_{2^{i + M} r + 2^i r}}{v_{2^i r}^2 v_{2^{i + M} r}^2} &= \sum_{i = 1}^M \frac{q^{2^i r} u_{2^{i + N} r - 2^i r} u_{2^{i + N} r + 2^i r}}{v_{2^i r}^2 v_{2^{i + N} r}^2}, \\
	\sum_{i = 1}^{2N} \frac{(\pm 1)^i q^{2^i r} u_{2^{i + 2M} r - 2^i r} u_{2^{i + 2M} r + 2^i r}}{u_{2^i r}^2 u_{2^{i + 2M} r}^2} &= \sum_{i = 1}^{2M} \frac{(\pm 1)^i q^{2^i r} u_{2^{i + 2N} r - 2^i r} u_{2^{i + 2N} r + 2^i r}}{u_{2^i r}^2 u_{2^{i + 2N} r}^2}, \\
	\sum_{i = 1}^{2N} \frac{(\pm 1)^i q^{2^i r} u_{2^{i + 2M} r - 2^i r} u_{2^{i + 2M} r + 2^i r}}{v_{2^i r}^2 v_{2^{i + 2M} r}^2} &= \sum_{i = 1}^{2M} \frac{(\pm 1)^i q^{2^i r} u_{2^{i + 2N} r - 2^i r} u_{2^{i + 2N} r + 2^i r}}{v_{2^i r}^2 v_{2^{i + 2N} r}^2}.
\end{align*}
\end{corollary}

\begin{corollary} \label{Cor: 32} Let $ r > 0 $ be an integer. Then
\begin{align*}
	\sum_{i = 1}^\infty \frac{4q^{2^i r} u_{2^i r} u_{3 \cdot 2^i r}}{u_{2^i r}^2 u_{2^{i + 1} r}^2} &= \frac{v_{2r}^2}{u_{2r}^2} - D, \\
	\sum_{i = 1}^\infty \frac{4q^{2^i r} u_{2^i r} u_{3 \cdot 2^i r}}{v_{2^i r}^2 v_{2^{i + 1} r}^2} &= \frac{1}{D} - \frac{u_{2r}^2}{v_{2r}^2}, \\
	\sum_{i = 1}^\infty \frac{(\pm 1)^i 4q^{2^i r} u_{3 \cdot 2^i r} u_{5 \cdot 2^i r}}{u_{2^i r}^2 u_{2^{i + 2} r}^2} &= \frac{v_{4r}^2}{u_{4r}^2} \pm \frac{v_{2r}^2}{u_{2r}^2} - D(1 \pm 1), \\
	\sum_{i = 1}^\infty \frac{(\pm 1)^i 4q^{2^i r} u_{3 \cdot 2^i r} u_{5 \cdot 2^i r}}{v_{2^i r}^2 v_{2^{i + 2} r}^2} &= -\frac{u_{4r}^2}{v_{4r}^2} \mp \frac{u_{2r}^2}{v_{2r}^2} + \frac{1}{D} (1 \pm 1).
\end{align*}
\end{corollary}

\begin{example} \label{Eg: 33} To apply to the case of Fibonacci and Lucas numbers, we set $ p = 1 $ and $ q = -1 $, then $ \alpha = \varphi $, $ \beta = 1 - \varphi $, and $ D = 5 $. Corollary \ref{Cor: 32} gives
\begin{align*}
	\sum_{i = 1}^\infty \frac{F_{3 \cdot 2^i r}}{F_{2^i r} F_{2^{i + 1} r}^2} &= \frac{L_{2r}^2}{4F_{2r}^2} - \frac{5}{4}, \\
	\sum_{i = 1}^\infty \frac{F_{2^i r} F_{3 \cdot 2^i r}}{L_{2^i r}^2 L_{2^{i + 1} r}^2} &= \frac{1}{20} - \frac{F_{2r}^2}{4L_{2r}^2}, \\
	\sum_{i = 1}^\infty \frac{(\pm 1)^i F_{3 \cdot 2^i r} F_{5 \cdot 2^i r}}{F_{2^i r}^2 F_{2^{i + 2} r}^2} &= \frac{L_{4r}^2}{4F_{4r}^2} \pm \frac{L_{2r}^2}{4F_{2r}^2} - \frac{5}{4} (1 \pm 1), \\
	\sum_{i = 1}^\infty \frac{(\pm 1)^i F_{3 \cdot 2^i r} F_{5 \cdot 2^i r}}{L_{2^i r}^2 L_{2^{i + 2} r}^2} &= -\frac{F_{4r}^2}{4L_{4r}^2} \mp \frac{F_{2r}^2}{4L_{2r}^2} + \frac{1}{20} (1 \pm 1). 
\end{align*}
\end{example}

\section{Further discussion} \label{Sec: 6}

It is possible to rewrite Theorems \ref{Th: 3} to \ref{Th: 6} in terms of $ \alpha $ or even consider the more general identity
\begin{align*}
	\frac{(A + B)(x_{i + j} y_i - x_i y_{i + j})}{(Ax_i - By_i)(Ax_{i + j} - By_{i + j})} = \frac{x_i + y_i}{Ax_i - By_i} - \frac{x_{i + j} + y_{i + j}}{Ax_{i + j} - By_{i + j}}.
\end{align*}
The results should be similar except that the limit does not vanish in this case. We leave it to the readers to work it out. 

Furthermore, we believe it is possible to further generalize the identity in \S \ref{Sec: 4} to include generalized Lucas sequences as well. This will be be an interesting topic to pursue in future research.

\end{document}